\documentclass[12pt]{amsart}
\usepackage{amssymb}
\usepackage{amsmath}
\usepackage{amscd}
\theoremstyle{plain}
\newtheorem{theorem}{Theorem}[section]
\newtheorem{corollary}[theorem]{Corollary}
\newtheorem{proposition}[theorem]{Proposition}
\newtheorem{lemma}[theorem]{Lemma}

{\theoremstyle{remark}

\newtheorem{remark}[theorem]{Remark}}
{\theoremstyle{definition}
\newtheorem{definition}[theorem]{Definition}

\newtheorem{example}[theorem]{Example}}

\newcommand{\rom}{\renewcommand{\labelenumi}{{\rm (\roman{enumi})}}%
\renewcommand{\itemsep}{0pt}}

\newcommand{\Z}{\mathbb{Z}}
\newcommand{\Q}{\mathbb{Q}}

\newcommand{\tM}{{\tilde{M}}}

\newcommand{\cC}{{\mathcal C}}

\DeclareMathOperator{\Hom}{Hom}
\DeclareMathOperator{\ab}{ab}
\DeclareMathOperator{\coker}{coker}

\begin{document}
\title[Permutation presentations of modules over finite groups]
{Permutation presentations of modules over finite groups.}
\author[Takeshi KATSURA]{Takeshi KATSURA}
\address{Department of Mathematics, 
Hokkaido University, Kita 10, Nishi 8, 
Kita-Ku, Sapporo, 060-0810, JAPAN}
\email{katsura@math.sci.hokudai.ac.jp}
\thanks{The author is partially supported by JSPS Research Fellow.} 
\date{}



\begin{abstract}
We introduce a notion of permutation presentations 
of modules over finite groups, 
and completely determine finite groups 
over which every module has a permutation presentation. 
To get this result, 
we prove that every coflasque module over a cyclic $p$-group 
is permutation projective. 
\end{abstract}

\maketitle

\section{Introduction}

In this paper, 
we introduce a notion of permutation presentations 
of modules over finite groups, 
and completely determine finite groups 
over which every module has a permutation presentation. 
This result is used in \cite{Ka} to solve a problem 
on $C^*$-algebras. 
Throughout this paper, 
$\Gamma$ denotes a finite group. 
Every modules are left modules.

\begin{definition}
A $\Z\Gamma$-module $F$ 
is called a {\it permutation module} 
if we can find a basis of $F$ as a $\Z$-module 
which is permuted by the action of $\Gamma$. 
\end{definition}

For a $\Gamma$-set $X$, 
the free abelian group $\Z[X]$ 
whose basis is given by $\{[x]\}_{x\in X}$ 
is a permutation $\Z\Gamma$-module, 
and every permutation module 
is isomorphic to this type of module. 

\begin{definition}
A {\em permutation presentation} of a $\Z\Gamma$-module $M$ is 
an exact sequence 
\[
\begin{CD}
0 @>>> F @>>> F @>>> M @>>> 0
\end{CD}
\]
of $\Z\Gamma$-modules 
where $F$ is a permutation $\Z\Gamma$-module. 
A permutation presentation is said to be {\em countable} 
when $F$ in the sequence above is countable. 
\end{definition}

If a $\Z\Gamma$-module $M$ 
has a countable permutation presentation, 
then $M$ is countable. 
The following proposition, 
which will be proved in Section~\ref{Sec:PP}, 
says that the converse also holds. 

\begin{proposition}\label{Prop:cou}
If a countable $\Z\Gamma$-module has a permutation presentation, 
then it has a countable permutation presentation. 
\end{proposition}

For a $\Gamma$-set $X$, 
an endomorphism $\varphi$ of the permutation $\Z\Gamma$-module $\Z[X]$ 
defines a $X\times X$-matrix $A=(A_{x,x'})_{x,x'\in X}$ 
with $A_{x,x'}\in\Z$ 
by $\varphi([x])=\sum_{x'\in X}A_{x,x'}[x']$. 
This matrix $A$ satisfies 
that $\{x'\in X\mid A_{x,x'}\neq 0\}$ is finite for all $x$, 
and $A_{x,x'}=A_{\gamma(x),\gamma(x')}$ for all $x,x'$ 
and $\gamma\in\Gamma$. 
Conversely such a $\Gamma$-invariant $X\times X$-matrix $A$ 
with integer entries defines an endomorphism of $\Z[X]$. 
Therefore a permutation presentation of a $\Z\Gamma$-module $M$ 
is presented by a $\Gamma$-set $X$ 
and a $\Gamma$-invariant $X\times X$-matrix $A$ with integer entries 
such that $\ker A=0$ and $\coker A\cong M$ as $\Z\Gamma$-modules. 
We will use this observation in \cite{Ka}. 

The following is the main theorem of this paper. 

\begin{theorem}\label{Thm:Main}
For a finite group $\Gamma$, 
every $\Z\Gamma$-module has a permutation presentation 
if and only if every Sylow subgroup of $\Gamma$ is cyclic. 
\end{theorem}

There exists an explicit description of finite groups 
all of whose Sylow subgroups are cyclic. 
Such a group is isomorphic to 
a semi-direct product $\Gamma=(\Z/m\Z)\rtimes (\Z/n\Z)$ 
such that $m$ and $n$ are relatively prime integers 
and that $\Z/m\Z$ is the commutator group of $\Gamma$ 
(for detail, see \cite[10.1.10]{Ro} for example). 
By Theorem~\ref{Thm:Main}, 
we see that 
for a finite cyclic group $\Gamma$ 
every $\Z\Gamma$-module has a permutation presentation. 
We remark that when the order of $\Gamma$ is prime, 
this fact follows from \cite[Theorem 1.3]{Sp}. 

The condition that a given $\Z\Gamma$-module $M$ has 
a permutation presentation 
is much weaker than it looks. 
This is equivalent to the condition that 
the permutation projective dimension of $M$ is $0$ or $1$, 
which was introduced by Arnold in \cite{A} 
(see the beginning of Section~\ref{Sec:PP} in this paper). 
In \cite{A2}, 
a finite group $\Gamma$ was said to have 
a global permutation projective dimension $1$ 
if every {\em finitely generated} $\Z\Gamma$-module $M$ 
has a permutation projective dimension $0$ or $1$. 
Using the results in \cite{EM}, 
Arnold showed that 
the condition that every Sylow subgroup is cyclic 
characterizes finite groups 
with global permutation projective dimension $1$ 
in his sense 
(see the remark after \cite[Definition 1.5]{A2}). 
On the other hand, 
one can define a global permutation projective dimension 
of a group $\Gamma$ using all $\Z\Gamma$-modules, 
and Theorem~\ref{Thm:Main} with the remark above 
shows that 
the condition for a finite group 
to have a global permutation projective dimension $1$
in this sense 
is the same as the one considered by Arnold. 
Thus two notions of global permutation projective dimension $1$ 
coincide. 
The author does not know whether this coincidence happens 
for higher dimensions
(see \cite{A} for the precise definition of 
permutation projective dimension). 

In order to prove Theorem \ref{Thm:Main}, 
we can use many arguments in \cite{A} and \cite{EM} 
without big changes. 
However, in the one step in \cite{EM}, 
Endo and Miyata seemed to use 
a property of 
modules over the Dedekind domain $\Z[\zeta_{q}]$ 
which is valid only for {\em finitely generated} modules, 
and hence we cannot use their argument. 
We need a new idea to complete this step 
in our case. 
To express this step, 
we introduce the following notions 
(see \cite[Subsection 2.10]{L}). 

\begin{definition}
A $\Z\Gamma$-module $M$ 
is said to be {\it permutation projective} 
if it is isomorphic to a direct summand of some permutation module, 
and {\it coflasque} if $M$ is free as a $\Z$-module 
and $H^1(\Gamma',M)=0$ for all subgroups $\Gamma'$ of $\Gamma$. 
\end{definition}

For a definition and results 
of the group cohomologies $H^i(-,-)$, 
we consult the books of Brown \cite{Br} 
and Serre \cite{Se}. 
By Shapiro's lemma, 
every permutation projective $\Z\Gamma$-module is coflasque 
(see Remark~\ref{Rem:Shapiro}). 
Now the following theorem is the most difficult part 
in our proof of Theorem~\ref{Thm:Main}. 

\begin{theorem}\label{Thm:TecMain}
Let $\Gamma$ be a cyclic $p$-group. 
Then every coflasque $\Z\Gamma$-module is 
permutation projective. 
\end{theorem}

This theorem seems to have its own importance. 
Recently, several authors investigate the structures 
of $\Z\Gamma$-modules which is free as a $\Z$-module. 
Such a module is called a {\em generalized lattice} in \cite{BCK}. 
The proof of Theorem \ref{Thm:TecMain} is inspired from \cite{BCK}, 
and is more elementary (though longer) than 
the arguments in \cite{EM}. 

This paper is organized as follows. 
In Section~\ref{Sec:Obs}, 
we give a cohomological obstruction for a module 
to have permutation presentations. 
We also prove 
the ``only if'' part of Theorem \ref{Thm:Main}. 
In Section~\ref{Sec:PP}, 
we sketch the proofs of propositions 
which were proved in \cite{A} and \cite{EM} 
for finitely generated modules 
explicitly or implicitly. 
These propositions give the induction 
of Theorem~\ref{Thm:Main} to Theorem~\ref{Thm:TecMain}. 
We also prove Proposition~\ref{Prop:cou}. 
Section~\ref{Sec:ProjMod}
is a preparatory section for the proof of Theorem~\ref{Thm:TecMain}, 
and in Section~\ref{Sec:CofMod} 
we prove Theorem~\ref{Thm:TecMain} 
and thus complete the proof of Theorem~\ref{Thm:Main}. 

\medskip

\noindent
{\bfseries Acknowledgments.} 
In March 2006, 
the author held a seminar about this paper 
at University of Tokyo, 
and he is grateful to the participants 
for asking stimulating questions and 
simplifying some arguments in an earlier draft of this paper.

\section{An obstruction to have permutation presentations}
\label{Sec:Obs}

In this section, 
we give a cohomological obstruction for a $\Z\Gamma$-module 
to have permutation presentations, 
and prove the ``only if'' part of Theorem~\ref{Thm:Main}. 

Recall that the exponent $\exp(\Gamma)$ of a finite group $\Gamma$ 
is the smallest positive integer $n$ satisfying 
$\gamma^{n}=1$ for all $\gamma\in \Gamma$. 
The following lemma is easy to see, 
hence we omit the proof. 

\begin{lemma}\label{Lem:nG}
The exponent $\exp(\Gamma)$ of $\Gamma$ divides the order $|\Gamma|$, 
and we have $\exp(\Gamma)=|\Gamma|$ if and only if 
every Sylow subgroup of $\Gamma$ is cyclic. 
\end{lemma}

\begin{lemma}\label{Lem:H2P}
For a permutation $\Z\Gamma$-module $F$, 
the group $H^2(\Gamma,F)$ is annihilated by $\exp(\Gamma)$. 
\end{lemma}

\begin{proof}
By the definition of permutation modules, 
we see that a permutation $\Z\Gamma$-module $F$ is isomorphic to 
a direct sum of modules in the form $\Z[\Gamma/\Gamma']$ 
for a subgroup $\Gamma'$ of $\Gamma$. 
For a finite group $\Gamma$, 
$H^i(\Gamma,-)$ commutes with taking direct sums. 
Hence $H^2(\Gamma,F)$ is isomorphic to 
a direct sum of $\{H^2(\Gamma,\Z[\Gamma/\Gamma'])\}$ 
for subgroups $\{\Gamma'\}$ of $\Gamma$. 
Using Shapiro's Lemma (cf.\ \cite[III.6.4]{Br}) 
and some well-known facts (cf.\ \cite[Exercise IV.3.3]{Br}), 
we obtain 
\[
H^2(\Gamma,\Z[\Gamma/\Gamma'])\cong 
H^2(\Gamma',\Z)\cong 
H^1(\Gamma',\Q/\Z)\cong \Hom(\Gamma',\Q/\Z)
\cong \widehat{(\Gamma')^{\ab}}
\cong (\Gamma')^{\ab}
\]
where $(\Gamma')^{\ab}$ is the abelianization of $\Gamma'$ 
and $\widehat{(\Gamma')^{\ab}}$ is its dual group. 
It is clear 
that an order of any element of $(\Gamma')^{\ab}$ divides $\exp(\Gamma)$. 
Hence $H^2(\Gamma,\Z[\Gamma/\Gamma'])$ 
is annihilated by $\exp(\Gamma)$ 
for any subgroup $\Gamma'$ of $\Gamma$. 
This shows that 
$H^2(\Gamma,F)$ is annihilated by $\exp(\Gamma)$. 
\end{proof}

\begin{remark}\label{Rem:Shapiro}
In a similar way to the proof of Lemma~\ref{Lem:H2P}, 
we can show $H^1(\Gamma,F)=0$ for a permutation $\Z\Gamma$-module $F$ 
because $H^1(\Gamma,\Z[\Gamma/\Gamma'])\cong 
H^1(\Gamma',\Z)\cong \Hom(\Gamma',\Z)=0$ 
for a subgroup $\Gamma'$ of $\Gamma$. 
From this fact, we see that 
every permutation projective $\Z\Gamma$-module is coflasque. 
\end{remark}

\begin{remark}
If $P$ is a permutation projective $\Z\Gamma$-module, 
then $P$ is coflasque and satisfies that 
$H^2(\Gamma',P)$ is annihilated by $\exp(\Gamma')$ 
for all subgroups $\Gamma'$ of $\Gamma$ by Lemma~\ref{Lem:H2P}. 
The author does not know 
whether the converse holds or not. 
This problem is closely related to the problem of 
determining whether a given module has a permutation presentation 
(see Corollary~\ref{Cor:NM}). 
\end{remark}

The following gives a cohomological obstruction for a $\Z\Gamma$-module 
to have permutation presentations. 

\begin{proposition}\label{Prop:obstr}
If a $\Z\Gamma$-module $M$ has a permutation presentation, 
then $H^1(\Gamma,M)$ is annihilated by $\exp(\Gamma)$. 
\end{proposition}

\begin{proof}
Let $0\to F\to F\to M\to 0$ 
be a permutation presentation of $M$.
Since $H^1(\Gamma,F)=0$, 
the long exact sequence of cohomologies 
gives us an injection $H^1(\Gamma,M)\to H^2(\Gamma,F)$. 
By Lemma \ref{Lem:H2P}, 
$H^2(\Gamma,F)$ is annihilated by $\exp(\Gamma)$. 
Hence so is $H^1(\Gamma,M)$. 
We are done. 
\end{proof}

By this proposition, 
we can see that if $H^1(\Gamma,M)$ is not annihilated by $\exp(\Gamma)$ 
then $M$ has no permutation presentations. 

\begin{example}\label{Ex1}
Let $\Gamma\cong (\Z/2\Z)\times (\Z/2\Z)$ be a group generated 
by two elements $\sigma,\tau\in\Gamma$ 
with relations $\sigma^2=\tau^2=1$ and $\sigma\tau=\tau\sigma$. 
We have $\exp(\Gamma)=2$. 
Let $M=(\Z/4\Z)^3$ be a $\Z\Gamma$-module 
where an action of $\Gamma$ on $M$ is defined by 
\[
\sigma\big((a,b,c)\big)=(b, a, -a-b-c),\quad 
\tau\big((a,b,c)\big)=(c,-a-b-c, a), 
\]
for $a,b,c\in \Z/4\Z$. 
One can compute $H^1(\Gamma,M)\cong \Z/4\Z$
which is not annihilated by $\exp(\Gamma)=2$. 
Hence by Proposition~\ref{Prop:obstr}, 
we see that $M$ has no permutation presentations. 
\end{example}

The following example shows that 
the obstruction given by Proposition~\ref{Prop:obstr}
is not complete. 

\begin{example}
Let $\Gamma\cong (\Z/2\Z)\times (\Z/2\Z)$ be as in Example~\ref{Ex1}. 
Let $M=\Z/8\Z$ be a $\Z\Gamma$-module 
where an action of $\Gamma$ on $M$ is defined by 
$\sigma(m)=-m$ and $\tau(m)=3m$ for $m\in M$. 
A straightforward computation shows $H^1(\Gamma,M)\cong\Z/2\Z$. 
Since $H^1(\Gamma,M)$ is annihilated by $\exp(\Gamma)=2$, 
we cannot conclude from Proposition~\ref{Prop:obstr} 
that $M$ has no permutation presentations. 
However, the following argument shows that 
$M$ has no permutation presentations. 
Let $F=\Z[M]$ be a permutation $\Z\Gamma$-module, 
and define $\pi\colon F\to M$ by $\pi([m])=m$. 
Let $N_M$ be the kernel of the surjection $\pi$. 
Then we will see in Corollary~\ref{Cor:NM} that 
$M$ has a permutation presentation 
if and only if 
the $\Z\Gamma$-module $N_M$ is 
permutation projective. 
We can show that $2H^2(\Gamma,N_M)$ coincides with 
the image of the injection $H^1(\Gamma,M)\to H^2(\Gamma,N_M)$ 
in the long exact sequence, 
and hence $2H^2(\Gamma,N_M)\neq 0$. 
Thus $H^2(\Gamma,N_M)$ is not annihilated by $\exp(\Gamma)=2$, 
and hence $N_M$ is not permutation projective 
by Lemma~\ref{Lem:H2P}. 
Therefore $M$ has no permutation presentations. 
\end{example}

We use Proposition~\ref{Prop:obstr} 
to show the ``only if'' part of Theorem~\ref{Thm:Main}. 
The proof is inspired by \cite[Theorem 1.5]{EM}. 

Recall that the augmentation ideal $I_\Gamma$ 
is the kernel of the surjection $\pi\colon \Z\Gamma\to\Z$ 
defined by $\pi(\gamma)=1$ for all $\gamma\in \Gamma$. 
Since $I_\Gamma$ is an ideal of $\Z\Gamma$, 
it is a $\Z\Gamma$-module. 
We can easily compute $H^1(\Gamma,I_\Gamma)\cong \Z/|\Gamma|\Z$. 

\begin{proposition}\label{Prop:Ing3}
If $I_\Gamma$ has a permutation presentation, 
then every Sylow subgroup of $\Gamma$ is cyclic. 
\end{proposition}

\begin{proof}
If $I_\Gamma$ has a permutation presentation, 
then $H^1(\Gamma,I_\Gamma)\cong \Z/|\Gamma|\Z$ 
is annihilated by $\exp(\Gamma)$ by Proposition \ref{Prop:obstr}. 
This shows that $\exp(\Gamma)$ is a multiple of $|\Gamma|$. 
By Lemma~\ref{Lem:nG}, 
this implies that every Sylow subgroup of $\Gamma$ is cyclic. 
\end{proof}

This proposition proves the ``only if'' part of Theorem~\ref{Thm:Main}.

\section{Coflasque modules and permutation projective modules}
\label{Sec:PP}

In this section, 
we see that we can reduce Theorem~\ref{Thm:Main} 
to Theorem~\ref{Thm:TecMain}. 
Most of the steps of the reduction 
were shown for finitely generated modules 
in \cite{A} and \cite{EM}, 
and their proofs can be adapted to our case 
without big changes. 
We review their arguments for readers' convenience 
and for the references in the rest of this paper. 
At the last of this section, 
we prove Proposition~\ref{Prop:cou}. 

In \cite{A}, 
Arnold introduced permutation projective dimensions 
(or ${\mathcal C}_1$ dimensions) 
of finitely generated $\Z\Gamma$-modules 
for a finite group $\Gamma$. 
His definition can be applied for an arbitrary $\Z\Gamma$-modules, 
and we can see that 
a $\Z\Gamma$-module $M$ has 
a permutation projective dimension $0$ or $1$ 
if and only if there exists a short exact sequence 
\[
\begin{CD}
0 @>>> P_1 @>>> P_0 @>>> M @>>> 0
\end{CD}
\]
such that $P_0$ and $P_1$ are permutation projective $\Z\Gamma$-modules 
(see \cite{A} for the precise definition). 
Note that we can take $P_0$ and $P_1$ to be finitely generated 
when $M$ is a finitely generated $\Z\Gamma$-module 
whose permutation projective dimension is $0$ or $1$ 
(see Proposition~\ref{Prop:ppd01}). 
Our first observation for the proof of Theorem~\ref{Thm:Main} 
is that this condition is equivalent to have a permutation presentation. 

\begin{proposition}\label{Prop:ppd}
A $\Z\Gamma$-module $M$ has a permutation presentation 
if and only if there exists a short exact sequence 
\[
\begin{CD}
0 @>>> P_1 @>>> P_0 @>>> M @>>> 0
\end{CD}
\]
such that $P_0$ and $P_1$ are permutation projective $\Z\Gamma$-modules. 
\end{proposition}

\begin{proof}
This follows from the next lemma. 
\end{proof}

\begin{lemma}\label{Lem:Swindle}
For two permutation projective $\Z\Gamma$-modules $P_0$ and $P_1$, 
there exists a permutation $\Z\Gamma$-module $F$ 
such that $P_0\oplus F\cong P_1\oplus F\cong F$. 
\end{lemma}

\begin{proof}
We use the technique called ``Eilenberg trick'' 
(cf.\ \cite[Lemma VIII.2.7]{Br}). 
Let $F_0$ and $F_1$ be permutation $\Z\Gamma$-modules 
containing $P_0$ and $P_1$ as direct summands, respectively. 
Set $F=\bigoplus_{k=1}^\infty(F_0\oplus F_1)$. 
This $\Z\Gamma$-module $F$ is a permutation module
satisfying $P_0\oplus F\cong P_1\oplus F\cong F$. 
\end{proof}

Recall that for a $\Z\Gamma$-module $M$ 
and a subgroup $\Gamma'$ of $\Gamma$, 
we have 
\[
H^0(\Gamma',M)\cong M^{\Gamma'}
:=\{m\in M\mid \text{$\gamma m=m$ for all $\gamma\in \Gamma'$}\}. 
\]

\begin{definition}
A $\Z\Gamma$-module map $\pi\colon M\to N$ 
between $\Z\Gamma$-modules $M,N$ 
is said to be {\em $H^0$-surjective} 
if its restriction $\pi\colon M^{\Gamma'}\to N^{\Gamma'}$ 
is surjective for every subgroup $\Gamma'$ of $\Gamma$. 
\end{definition}

By setting $\Gamma'=\{1\}$, 
we see that an $H^0$-surjective map is surjective. 

\begin{lemma}\label{Lem:123}
Let 
\[
\begin{CD}
0 @>>> M_1 @>>> M_2 @>\pi>> M_3 @>>> 0.
\end{CD}
\]
be an exact sequence of $\Z\Gamma$-modules. 
Then the following hold. 
\begin{enumerate}
\item If $M_1$ is coflasque, 
then $\pi$ is $H^0$-surjective. 
\item 
If $\pi$ is $H^0$-surjective 
and $M_2$ is coflasque, 
then $M_1$ is coflasque. 
\item 
If $\pi$ is $H^0$-surjective 
and $M_3$ is permutation projective, 
then this sequence splits. 
\item If $M_1$ is coflasque and $M_3$ is permutation projective, 
then this sequence splits. 
\end{enumerate}
\end{lemma}

\begin{proof}
\begin{enumerate}
\item 
This follows from the long exact sequence of cohomologies. 
\item 
This follows from the long exact sequence of cohomologies, 
and the fact that a submodule of a free $\Z$-module is free. 
\item 
We may assume that $M_3$ is a permutation module, 
and in this case the conclusion easily follows from 
the definition of $H^0$-surjective maps. 
\item 
This follows from (1) and (3). 
\end{enumerate}
\end{proof}

\begin{proposition}\label{Prop:extpp}
Let 
\[
\begin{CD}
0 @>>> M_1 @>>> M_2 @>>> M_3 @>>> 0.
\end{CD}
\]
be an exact sequence of $\Z\Gamma$-modules. 
If both $M_1$ and $M_3$ are permutation projective, 
then the sequence splits, 
and hence $M_2\cong M_1\oplus M_3$ 
is also permutation projective. 
\end{proposition}

\begin{proof}
This follows from Lemma~\ref{Lem:123}~(4) 
because permutation projective modules are coflasque. 
\end{proof}

\begin{proposition}[{cf.\ \cite[Theorem 3.5]{A}}]\label{Prop:ppd01}
A $\Z\Gamma$-module $M$ 
has a permutation presentation 
if and only if 
for every permutation projective $\Z\Gamma$-module $P$ 
and every $H^0$-surjective map $\pi\colon P\to M$, 
$\ker\pi$ is permutation projective. 
\end{proposition}

\begin{proof}
The proof of the ``only if'' part goes 
similarly to the proof of ``Schanuel's lemma'' 
(\cite[Lemma VIII.4.2]{Br}). 
Suppose that a $\Z\Gamma$-module $M$ 
has a permutation presentation 
\[
\begin{CD}
0 @>>> F @>>> F @>\varphi>> M @>>> 0
\end{CD}
\]
where $F$ is a permutation $\Z\Gamma$-module. 
Take a permutation projective $\Z\Gamma$-module $P$ 
and an $H^0$-surjective map $\pi\colon P\to M$. 
Let 
\[
N:=\{(f,p)\in F\times P\mid \varphi(f)=\pi(p)\}
\]
be the pull-buck of the two surjections $\varphi$ and $\pi$. 
Then we have the following commutative diagram 
with exact rows and columns: 
\[
\begin{CD}
@. @. 0 @. 0\\
@. @. @VVV @VVV @.\\
@. @. \ker \pi @= \ker\pi \\
@. @. @VVV @VVV @.\\
0 @>>> F @>>> N @>>> P @>>> 0\\
@. @| @VVV @V\pi VV @.\\
0 @>>> F @>>> F @>\varphi>> M @>>> 0\\ 
@. @. @VVV @VVV @.\\
@. @. 0 @. 0
\end{CD}
\]
Since both $F$ and $P$ are permutation projective, 
the middle row splits and $N\cong F\oplus P$ 
is permutation projective 
by Proposition~\ref{Prop:extpp}. 
Since $\pi$ is $H^0$-surjective 
and $P$ is permutation projective, 
$\ker\pi$ is coflasque by Lemma~\ref{Lem:123}~(2). 
By Lemma~\ref{Lem:123}~(4) 
the middle column splits 
and we get $N\cong \ker\pi \oplus F$. 
This shows that $\ker\pi$ is permutation projective. 

Conversely suppose that 
$\ker\pi$ is permutation projective 
for every permutation projective $\Z\Gamma$-module $P$ 
and every $H^0$-surjective map $\pi\colon P\to M$. 
Let us define $\pi\colon \Z[M]\to M$ by $\pi([m])=m$. 
Since $\Z[M]$ is a permutation $\Z\Gamma$-module 
and $\pi$ is $H^0$-surjective, 
$\ker\pi$ is permutation projective. 
Thus by Proposition~\ref{Prop:ppd} 
$M$ has a permutation presentation. 
\end{proof}

For a $\Z\Gamma$-module $M$, 
let $N_M$ be the kernel of 
the $H^0$-surjective map $\pi\colon \Z[M]\to M$ 
defined by $\pi([m])=m$ 
which is considered in the proof of 
Proposition~\ref{Prop:ppd01}. 
Then by Lemma~\ref{Lem:123}~(2), 
$N_M$ is coflasque. 
By Proposition~\ref{Prop:ppd01}, 
we get the following. 

\begin{corollary}\label{Cor:NM}
A $\Z\Gamma$-module $M$ 
has a permutation presentation 
if and only if 
the coflasque $\Z\Gamma$-module $N_M$ is 
permutation projective. 
\end{corollary}

\begin{lemma}\label{lem:embed}
For a $\Z\Gamma$-module $N$ which is free as a $\Z$-module, 
there exist a permutation $\Z\Gamma$-module $F$ 
and an injective $\Z\Gamma$-module map $\iota\colon N\to F$. 
\end{lemma}

\begin{proof}
Consider $N$ as a $\Z$-module, 
and set a tensor product $F=(\Z\Gamma)\otimes N$ as a $\Z$-module. 
The natural action of $\Gamma$ on $\Z\Gamma$ 
defines an action of $\Gamma$ on $F$. 
By this action, $F$ is a free $\Z\Gamma$-module. 
The map $\iota\colon N\to F$ defined by 
$\iota(n)=\sum_{\gamma\in \Gamma}[\gamma]\otimes (\gamma^{-1}n)$ 
is an injective $\Z\Gamma$-module map. 
\end{proof}

\begin{proposition}\label{Prop:Ing1}
Every $\Z\Gamma$-module 
has a permutation presentation 
if and only if 
every coflasque $\Z\Gamma$-module is permutation projective. 
\end{proposition}

\begin{proof}
Suppose that 
every $\Z\Gamma$-module 
has a permutation presentation. 
Take a coflasque $\Z\Gamma$-module $N$. 
By Lemma~\ref{lem:embed}, 
there exist a permutation $\Z\Gamma$-module $F$ 
and an injective $\Z\Gamma$-module map $\iota\colon N\to F$. 
Set a $\Z\Gamma$-module $M$ by $M=F/\iota(N)$, 
and denote by $\pi$ the natural surjection $F\to M$. 
Since $N$ is coflasque, 
$\pi$ is $H^0$-surjective by Lemma~\ref{Lem:123}~(1). 
Then by Proposition~\ref{Prop:ppd01}, 
$\ker\pi\cong N$ is permutation projective. 
This proves the ``only if'' part. 
The ``if'' part follows from Corollary~\ref{Cor:NM}. 
\end{proof}

\begin{lemma}[{\cite[Lemma 1.4]{EM}}]\label{Lem:Sylow}
A $\Z\Gamma$-module $M$ is permutation projective 
if and only if 
$M$ is a permutation projective 
$\Z\Gamma'$-module for every Sylow subgroup $\Gamma'$ of $\Gamma$. 
\end{lemma}

\begin{proof}
Clearly, 
if $M$ is a permutation projective $\Z\Gamma$-module, 
then it is a permutation projective 
$\Z\Gamma'$-module for every subgroup $\Gamma'$ of $\Gamma$. 
Let $M$ be a $\Z\Gamma$-module 
which is permutation projective 
as a $\Z\Gamma'$-module for every Sylow subgroup $\Gamma'$ of $\Gamma$. 
Let $\pi\colon \Z[M]\to M$ 
be an $H^0$-surjective map defined by $\pi([m])=m$. 
By Lemma~\ref{Lem:123}~(3), 
the surjection $\pi\colon \Z[M]\to M$ 
has a $\Z\Gamma'$-module splitting map 
for every Sylow subgroup $\Gamma'$ of $\Gamma$. 
Then it has a $\Z\Gamma$-module splitting map 
(cf. \cite[Proposition III.2.2 and Theorem III.10.3]{Br}). 
Thus $M$ is permutation projective. 
\end{proof}

The analogous statement of the lemma above 
for coflasque modules also holds. 

\begin{proposition}\label{Prop:Ing2}
For a finite group $\Gamma$, 
every coflasque $\Z\Gamma$-module is permutation projective 
if and only if 
every coflasque $\Z\Gamma'$-module is permutation projective 
for every Sylow subgroup $\Gamma'$ of $\Gamma$. 
\end{proposition}

\begin{proof}
Suppose that 
every coflasque $\Z\Gamma$-module is permutation projective. 
Take a subgroup $\Gamma'$ of $\Gamma$ 
and a coflasque $\Z\Gamma'$-module $N$. 
Let $M=\Z\Gamma\otimes_{\Z\Gamma'}N$ be 
the induced $\Z\Gamma$-module from $N$. 
Using Shapiro's lemma, 
one can check that $M$ is a coflasque $\Z\Gamma$-module. 
Hence by the assumption, 
$M$ is a permutation projective $\Z\Gamma$-module. 
Considering as a $\Z\Gamma'$-module, 
$M$ is also permutation projective 
and contains $N$ as a direct summand. 
Hence $N$ is a permutation projective $\Z\Gamma'$-module. 
This proves the ``only if'' part. 

Conversely suppose that 
every coflasque $\Z\Gamma'$-module is permutation projective 
for every Sylow subgroup $\Gamma'$ of $\Gamma$. 
Take a coflasque $\Z\Gamma$-module $M$. 
For every Sylow subgroup $\Gamma'$ of $\Gamma$, 
$M$ is a permutation projective $\Z\Gamma'$-module 
because $M$ is coflasque as a $\Z\Gamma'$-module. 
Hence by Lemma~\ref{Lem:Sylow}, 
$M$ is a permutation projective $\Z\Gamma$-module. 
This proves the ``if'' part. 
\end{proof}

We have already seen the ``only if'' part of Theorem~\ref{Thm:Main} 
in Proposition~\ref{Prop:Ing3}, 
and the ``if'' part follows from Theorem~\ref{Thm:TecMain} 
by Proposition~\ref{Prop:Ing1} 
and Proposition~\ref{Prop:Ing2}. 
We will prove Theorem~\ref{Thm:TecMain} 
in Section~\ref{Sec:CofMod}. 
We finish this section by proving Proposition~\ref{Prop:cou}. 

\begin{lemma}\label{Lem:cou}
A countable permutation projective module $N$ is 
a direct summand of a countable permutation module. 
\end{lemma}

\begin{proof}
Let $N$ be a countable permutation projective $\Z\Gamma$-module. 
Choose a permutation $\Z\Gamma$-module $\Z[X]$ 
which contains $N$ as a direct summand, 
where $X$ is a $\Gamma$-set. 
For each $n\in N$, 
let us define a finite subset $Y_n\subset X$ by 
$Y_n:=\{x\in X\mid n_x\neq 0\}$ 
where $n=\sum_{x\in X}n_x[x]$. 
Then the set $Y=\bigcup_{n\in N}Y_n$ is a countable subset of $X$ 
which is globally invariant under the action of $\Gamma$. 
Set $F:=\Z[Y]\subset \Z[X]$. 
Then $F$ is a countable permutation module 
containing $N$ 
as a direct summand. 
\end{proof}

\begin{proof}[Proof of Proposition~\ref{Prop:cou}]
Let $M$ be a countable $\Z\Gamma$-module 
which has a permutation presentation. 
Then by Corollary~\ref{Cor:NM}, 
the coflasque $\Z\Gamma$-module $N_M$, 
which is defined to be the kernel of 
the $H^0$-surjective map $\pi\colon \Z[M]\to M$, 
is permutation projective. 
Since $\Z[M]$ and its submodule $N_M$ are countable, 
one can take a countable permutation module $F$ 
such that $\Z[M]\oplus F\cong N_M\oplus F\cong F$ 
by Lemma~\ref{Lem:cou} and the proof of Lemma~\ref{Lem:Swindle}. 
Thus we get a countable permutation presentation 
\[
\begin{CD}
0 @>>> F @>>> F @>>> M @>>> 0
\end{CD}
\]
of $M$. 
\end{proof}

\section{Projective modules over cyclic $p$-groups}\label{Sec:ProjMod}

Let $p$ be a prime number, 
and $q$ be a power of $p$. 
Let $\Gamma$ be a cyclic group of order $q$ 
whose generator is given by $\sigma$. 
We set $s\in\Z\Gamma$ by $s=\sum_{j=0}^{q-1}\sigma^j$. 
The group cohomology of a $\Z\Gamma$-module $M$ has a period $2$, i.e. 
$H^{k}(\Gamma,M)=H^{k+2}(\Gamma,M)$ for $k\geq 1$. 
Hence we are only interested in $H^1(\Gamma,M)$ and $H^2(\Gamma,M)$ 
which are easily computed as follows; 
\[
H^1(\Gamma,M)\cong \{m\in M\mid sm=0\}/(1-\sigma)M, \quad
H^2(\Gamma,M)\cong M^\sigma/sM,
\]
where $M^\sigma=\{m\in M\mid \sigma m=m\}$. 
The next proposition is well-known. 

\begin{proposition}\label{Prop:proj}
A $\Z\Gamma$-module $M$ is projective 
if and only if 
$M$ is free as a $\Z$-module and 
$H^1(\Gamma,M)=H^2(\Gamma,M)=0$. 
\end{proposition}

For a proof of this proposition, 
see \cite[Proposition 3.6]{Ri} 
where a much more general statement was proved. 

We use the following lemma 
in the next section. 

\begin{lemma}\label{Lem:FindF}
Let $M$ be a free $\Z\Gamma$-module. 
Suppose that 
a $\Z\Gamma$-submodule $N$ of $M$ satisfies that 
$p^lM\subset N$ and 
$H^1(\Gamma,p^{l-1}M\cap N)=0$ for some integer $l$. 
Then there exists a free $\Z\Gamma$-submodule $F\subset M$ 
satisfying that 
\begin{itemize}
\item $M/F$ is free as a $\Z$-module, 
\item $p^{l-1}M\subset F+N$, and 
\item $F\cap N=p^lF$. 
\end{itemize}
\end{lemma}

\begin{proof}
We set $M'=(p^{l-1}M)/(p^{l-1}M\cap N)$ and 
let $\pi\colon p^{l-1}M\to M'$ be the natural surjection: 
\[
\begin{CD}
0 @>>> p^{l-1}M\cap N @>>> p^{l-1}M @>\pi>> M' @>>> 0.
\end{CD}
\]
Since $p^lM\subset N$, 
$M'$ is a $(\Z/p\Z)\Gamma$-module. 
Since $H^2(\Gamma,p^{l-1}M)=0$ 
and $H^3(\Gamma,p^{l-1}M\cap N)=H^1(\Gamma,p^{l-1}M\cap N)=0$, 
we have $H^2(\Gamma,M')=0$. 
Let $X$ be a basis of $M$ as a $\Z\Gamma$-module. 
Then $s\pi(p^{l-1}X)$ generates $sM'$ 
which coincides with $(M')^{\sigma}$ 
because $H^2(\Gamma,M')=0$. 
Choose a subset $Y\subset X$ such that 
$s\pi(p^{l-1}Y)$ is a basis of $(M')^{\sigma}$ 
as a $\Z/p\Z$-module. 
Let $F:=(\Z\Gamma)Y$ be the free $\Z\Gamma$-submodule of $M$ 
generated by $Y$. 
Then $M/F$ is free as a $\Z$-module 
(actually as a $\Z\Gamma$-module). 
We set $Y_0=\bigcup_{k=0}^{q-1}(\sigma-1)^kY$ 
which is a basis of $F$ as a $\Z$-module. 
By noticing $s=(\sigma-1)^{q-1}$ in $(\Z/p\Z)\Gamma$, 
we see that 
$\pi(p^{l-1}Y_0)=\bigcup_{k=0}^{q-1}(\sigma-1)^k\pi(p^{l-1}Y)$ 
is a basis of $M'$ 
as a $\Z/p\Z$-module. 
This shows that 
$p^{l-1}M\subset p^{l-1}F+N \subset F+N$. 
This also shows that 
if $x\in p^{l-1}F$ satisfies $\pi(x)=0$ 
then $x\in p^{l}F$. 
The proof ends once we show $F\cap N=p^lF$. 
Clearly we have $F\cap N\supset p^lF$. 
To derive a contradiction, 
assume that there exists $x_0\in (F\cap N)\setminus p^lF$. 
Choose a natural number $k$ 
such that $p^kx_0\in p^{l-1}F\setminus p^lF$. 
Since $p^kx_0\in N$, we have $\pi(p^kx_0)=0$. 
This contradicts the fact above. 
Thus we have $F\cap N=p^lF$. 
We are done. 
\end{proof}

\section{Coflasque modules over cyclic $p$-groups}\label{Sec:CofMod}

In this section, 
we show Theorem~\ref{Thm:TecMain}, 
and complete the proof of Theorem~\ref{Thm:Main}. 

Let $p$ be a prime number, 
and $K$ be a positive integer. 
Let $\Gamma$ be the cyclic group of order $p^K$ 
whose generator is given by $\sigma\in \Gamma$. 
Let us take $k\in\{0,1,\ldots,K\}$. 
We set $\sigma_k=\sigma^{p^k}$. 
Thus $\sigma_0=\sigma$, $\sigma_K=1$ 
and $\sigma_{k+1}=(\sigma_k)^p$. 
Let $\Gamma_k\subset \Gamma$ be 
the subgroup generated by $\sigma_k$. 
We have $\Gamma_k\cong \Z/p^{K-k}\Z$ 
and $\Gamma/\Gamma_k\cong \Z/p^{k}\Z$. 
For $l\in\{k,\ldots, K\}$, 
we set an element $s_k^{(l)}\in \Z\Gamma$ 
by $s_k^{(l)}=\sum_{j=0}^{p^{l-k}-1}\sigma_k^j$. 
We have 
\[
(1-\sigma_k)s_k^{(l)}=1-\sigma_l,\quad 
s_k^{(l)}s_l^{(m)}=s_k^{(m)} 
\]
for $k\leq l\leq m$. 
When $l=K$, 
the element $s_k^{(K)}\in \Z\Gamma$ 
is simply denoted by $s_k$. 

Every $\Z\Gamma$-module $M$ can also 
be considered as a $\Z\Gamma_k$-module for $k=0,1,\ldots,K$, 
and we have 
\[
H^1(\Gamma_k,M)\cong \{m\in M\mid s_km=0\}/(1-\sigma_k)M,\quad
H^2(\Gamma_k,M)\cong M^{\sigma_k}/s_kM.
\]
By a $\Z(\Gamma/\Gamma_k)$-module, 
we mean 
a $\Z\Gamma$-module $M$ with $M^{\sigma_k}=M$. 
Since $\{\Gamma_k\}_{k=0}^K$ exhausts all subgroups of $\Gamma$, 
a permutation $\Z\Gamma$-module is nothing but 
a direct sum of free $\Z(\Gamma/\Gamma_k)$-modules 
for $k=0,1,\ldots,K$.

\begin{definition}
Let $k=0,1,\ldots,K$. 
We write $M\in \cC_k$ 
if $M$ is a coflasque $\Z\Gamma$-module 
such that 
$M^{\sigma_k}$ is a projective $\Z(\Gamma/\Gamma_k)$-module. 
\end{definition}

Note that $M\in \cC_0$ if and only if $M$ is a coflasque $\Z\Gamma$-module, 
and $M\in \cC_K$ if and only if $M$ is a projective $\Z\Gamma$-module. 

\begin{lemma}\label{Lem:Vk}
Let $k\in\{1,2,\ldots,K\}$. 
For a coflasque $\Z\Gamma$-module $M$, 
the following conditions are equivalent; 
\begin{enumerate}
\rom 
\item $M\in \cC_k$, 
\item $p^{K-k}M^{\sigma_{k-1}}\subset s_{k-1}M$, 
\item $M^{\sigma_{k-1}}=s_{k-1}^{(k)}M^{\sigma_{k}}$, 
\item $M^{\sigma_0}=s_0^{(k)}M^{\sigma_{k}}$. 
\end{enumerate}
\end{lemma}

\begin{proof}
(i)$\Rightarrow$(ii): 
If $M^{\sigma_k}$ is a projective $\Z(\Gamma/\Gamma_k)$-module, 
then we have $M^{\sigma_{k-1}}=s_{k-1}^{(k)}M^{\sigma_k}$. 
Hence we get 
\[
p^{K-k}M^{\sigma_{k-1}}
=p^{K-k}s_{k-1}^{(k)}M^{\sigma_k}
=s_{k-1}M^{\sigma_k}
\subset s_{k-1}M.
\]

(ii)$\Rightarrow$(iii): 
Take $m\in M^{\sigma_{k-1}}$. 
By (ii), there exists $x\in M$ 
with $p^{K-k}m=s_{k-1}x$. 
Then we get $s_k(m-s_{k-1}^{(k)}x)=0$. 
Since $H^1(\Gamma_{k},M)=0$, 
we can find $y\in M$ such that 
$m-s_{k-1}^{(k)}x=(1-\sigma_{k})y=s_{k-1}^{(k)}(1-\sigma_{k-1})y$. 
Set $m'=x+(1-\sigma_{k-1})y\in M$. 
Then we have $m=s_{k-1}^{(k)}m'$. 
Since $m\in M^{\sigma_{k-1}}$, we see $m'\in M^{\sigma_{k}}$. 
Thus we have shown 
$M^{\sigma_{k-1}}\subset s_{k-1}^{(k)}M^{\sigma_{k}}$. 
The converse inclusion is obvious. 

(iii)$\Rightarrow$(iv): 
First we show that $M^{\sigma_{l}}=s_{l}^{(l+1)}M^{\sigma_{l+1}}$ 
holds for $l=0,1,\ldots,k-1$. 
By (iii), this holds for $l=k-1$. 
Suppose that this holds for $l$, 
and we show the equality for $l-1$. 
Take $m\in M^{\sigma_{l-1}}$. 
By the assumption, 
there exists $x\in M^{\sigma_{l+1}}$ 
such that $m=s_{l}^{(l+1)}x$. 
Then 
\[
s_{l}^{(l+1)}(m-s_{l-1}^{(l)}x)
=pm-s_{l-1}^{(l)}m=0
\]
because $m\in M^{\sigma_{l-1}}$. 
Since $H^1(\Gamma_l,M)=0$, 
we can find $y\in M$ with 
$m-s_{l-1}^{(l)}x=(1-\sigma_l)y=s_{l-1}^{(l)}(1-\sigma_{l-1})y$. 
Set $m'=x+(1-\sigma_{l-1})y$. 
Then we have $m=s_{l-1}^{(l)}m'$. 
The fact $m\in M^{\sigma_{l-1}}$ implies $m'\in M^{\sigma_{l}}$. 
Thus we get $M^{\sigma_{l-1}}\subset s_{l-1}^{(l)}M^{\sigma_{l}}$, 
and the converse inclusion is obvious. 
We showed that $M^{\sigma_{l}}=s_{l}^{(l+1)}M^{\sigma_{l+1}}$ 
holds for $l=0,1,\ldots,k-1$. 
Hence 
\[
M^{\sigma_{0}}
=s_0^{(1)}M^{\sigma_{1}}
=s_0^{(1)}s_1^{(2)}M^{\sigma_{2}}
=\cdots\cdots
=s_0^{(1)}s_1^{(2)}\cdots s_{k-1}^{(k)}M^{\sigma_{k}}
=s_0^{(k)}M^{\sigma_{k}}.
\]

(iv)$\Rightarrow$(i): 
The condition (iv) is nothing but 
$H^2(\Gamma/\Gamma_{k},M^{\sigma_{k}})=0$. 
On the other hand, 
we have $H^1(\Gamma/\Gamma_{k},M^{\sigma_{k}})=0$ 
for any coflasque $\Z\Gamma$-module $M$ 
because there exists 
an injection 
$H^1(\Gamma/\Gamma_{k},M^{\sigma_{k}})\to H^1(\Gamma,M)$ 
induced by the inflation 
(see \cite[Proposition VII.4]{Se}). 
Thus $M^{\sigma_{k}}$ is a projective $\Z(\Gamma/\Gamma_k)$-module 
by Proposition~\ref{Prop:proj}. 
\end{proof}

\begin{lemma}\label{Lem:canapply}
Let $k\in\{0,1,\ldots,K-1\}$. 
Let $\tM\in \cC_k$, 
and set $M=\tM^{\sigma_k}$ and $N=s_k\tM$ 
which are $\Z(\Gamma/\Gamma_k)$-modules. 
Then $p^{K-k}M\subset N\subset M$
and $H^1(\Gamma/\Gamma_k,p^{K-k-1}M\cap N)=0$. 
\end{lemma}

\begin{proof}
Clearly we have $p^{K-k}M\subset N\subset M$. 
Take $x\in p^{K-k-1}M\cap N$ with $s_0^{(k)}x=0$. 
Since $M$ is a projective $\Z(\Gamma/\Gamma_k)$-module, 
we have $H^1(\Gamma_0/\Gamma_k,p^{K-k-1}M)=0$. 
Hence there exists $y\in M$ with 
$x=(1-\sigma_0)p^{K-k-1}y$.
Since $x\in N=s_k\tM$, 
there exists $x=s_km$ with $m\in \tM$. 
We have 
\[
s_0 m=s_0^{(k)}x=s_0^{(k)}(1-\sigma_0)p^{K-k-1}y
=(1-\sigma_k)p^{K-k-1}y=0
\]
because $y\in M=\tM^{\sigma_k}$. 
Since $H^1(\Gamma_0,\tM)=0$, 
there exists $z\in \tM$ 
with $m=(1-\sigma_0)z$. 
Then 
\[
s_{k+1}(1-\sigma_0)(y-s_k^{(k+1)}z)
=(1-\sigma_0)p^{K-k-1}y-s_{k}(1-\sigma_0)z
=x-x
=0. 
\]
Since $H^1(\Gamma_{k+1},\tM)=0$, 
there exists $m_0\in \tM$ such that 
\[
(1-\sigma_0)(y-s_k^{(k+1)}z)
=(1-\sigma_{k+1})m_0
=(1-\sigma_0)s_k^{(k+1)}s_0^{(k)}m_0.
\]
Set $z'=z+s_0^{(k)}m_0\in \tM$. 
Then we have $(1-\sigma_0)y=(1-\sigma_0)s_k^{(k+1)}z'$. 
Since $y\in M=\tM^{\sigma_k}$, 
$s_k^{(k+1)}z'\in \tM^{\sigma_k}$. 
Hence the element 
\[
x':=p^{K-k-1}s_k^{(k+1)}z'
=s_{k+1}(s_k^{(k+1)}z')=s_kz' 
\]
is in $p^{K-k-1}M\cap N$. 
We have 
\[
(1-\sigma_0)x'=(1-\sigma_0)p^{K-k-1}(s_k^{(k+1)}z')
=(1-\sigma_0)p^{K-k-1}y=x
\]
This shows $H^1(\Gamma/\Gamma_k,p^{K-k-1}M\cap N)=0$. 
We are done. 
\end{proof}

\begin{proposition}\label{Prop:ind}
Let $k\in\{0,1,\ldots,K-1\}$. 
For $M\in \cC_k$, 
there exists a short exact sequence of $\Z\Gamma$-modules 
\[
\begin{CD}
0 @>>> F @>>> M\oplus P @>>> M' @>>> 0
\end{CD}
\]
where $F$ is a free $\Z(\Gamma/\Gamma_k)$-module, 
$P$ is a projective $\Z(\Gamma/\Gamma_k)$-module, 
and $M'\in \cC_{k+1}$. 
\end{proposition}

\begin{proof}
Since $M\in \cC_k$, 
$M^{\sigma_k}$ is a projective $\Z(\Gamma/\Gamma_k)$-module. 
Hence there exists a projective $\Z(\Gamma/\Gamma_k)$-module $P$ 
such that $M^{\sigma_k}\oplus P$ is a free $\Z(\Gamma/\Gamma_k)$-module. 
Set $\tM:=M\oplus P$. 
Then $\tM$ is a coflasque $\Z\Gamma$-module, 
and $\tM^{\sigma_k}=M^{\sigma_k}\oplus P$ 
is a free $\Z(\Gamma/\Gamma_k)$-module. 
By Lemma~\ref{Lem:canapply}, 
we can apply Lemma~\ref{Lem:FindF} 
to get a free $\Z(\Gamma/\Gamma_k)$-module $F\subset \tM^{\sigma_k}$ 
such that $\tM^{\sigma_k}/F$ is free as a $\Z$-module, 
$p^{K-k-1}\tM^{\sigma_k}\subset F+s_k\tM$ 
and $F\cap s_k\tM=p^{K-k}F$. 
Set $M'=\tM/F$. 
The proof ends if we show $M'\in \cC_{k+1}$. 
Since $\tM^{\sigma_k}/F$ and 
$\tM/\tM^{\sigma_k}\cong (1-\sigma_k)\tM\subset \tM$ 
are free as $\Z$-modules, 
$M'=\tM/F$ is also free as a $\Z$-module. 
To show that $H^1(\Gamma_l,M')=0$ for $l=0,1,\ldots,K$, 
it suffices to see 
by the long exact sequence of cohomologies 
that the map $H^2(\Gamma_l,F)\to H^2(\Gamma_l,\tM)$ is injective. 
Let us take $f\in F^{\sigma_{l}}$ such that $f=s_lm$ for some $m\in \tM$. 
We will show $f\in s_lF$. 
First consider the case $l\leq k$. 
Then $f=s_lm=s_k(s_l^{(k)}m)\in s_k\tM$. 
From $F\cap s_k\tM=p^{K-k}F$ 
we can find $f'\in F$ with $f=p^{K-k}f'$. 
Since $f\in F^{\sigma_{l}}$, $f'$ is also in $F^{\sigma_{l}}$. 
Since $F$ is a free $\Z(\Gamma/\Gamma_k)$-module, 
we can find $g\in F$ such that $f'=s_l^{(k)}g$. 
Then $s_lg=s_kf'=p^{K-k}f'=f$. 
Thus we have shown $f\in s_lF$. 
Next we consider the case $l> k$. 
We have $p^{l-k}f=s_k^{(l)}f\in s_k\tM$. 
Hence there exists $g\in F$ with $p^{l-k}f=p^{K-k}g$. 
Then we get $f=p^{K-l}g=s_lg\in s_lF$. 
Thus we have shown that 
the map $H^2(\Gamma_l,F)\to H^2(\Gamma_l,\tM)$ is injective. 
Therefore $M'$ is a coflasque $\Z\Gamma$-module. 
Since $F$ is a free $\Z(\Gamma/\Gamma_k)$-module, 
we have $H^1(\Gamma_k,F)=0$. 
This implies that $\tM^{\sigma_k}\to (M')^{\sigma_k}$ is surjective. 
From this fact and $p^{K-k-1}\tM^{\sigma_k}\subset F+s_k\tM$, 
we get $p^{K-k-1}(M')^{\sigma_k}\subset s_kM'$. 
Therefore Lemma~\ref{Lem:Vk}
shows $M'\in \cC_{k+1}$. 
This completes the proof. 
\end{proof}

\begin{proof}[Proof of Theorem~\ref{Thm:TecMain}]
Take a coflasque $\Z\Gamma$-module $M_0$. 
We have $M_0\in \cC_{0}$. 
By applying Proposition~\ref{Prop:ind} inductively, 
we get $\Z\Gamma$-modules 
$\{F_k, P_k\}_{k=0}^{K-1}$ and $\{M_k\}_{k=1}^{K}$ 
where $F_k$ is a free $\Z(\Gamma/\Gamma_k)$-module, 
$P_k$ is a projective $\Z(\Gamma/\Gamma_k)$-module
and $M_{k}\in \cC_{k}$ 
such that there exist short exact sequences 
\renewcommand{\theequation}{\fnsymbol{equation}}
\begin{equation}\label{FMPM}
\begin{CD}
0 @>>> F_k @>>> M_k\oplus P_k @>>> M_{k+1} @>>> 0. 
\end{CD}
\end{equation}
Since $M_K\in \cC_K$ is a projective $\Z\Gamma$-module, 
there exists a projective $\Z\Gamma$-module $P_K$ 
such that $F_K=M_K\oplus P_K$ is a free $\Z\Gamma$-module. 
We will show that 
\[
M_k\oplus \bigoplus_{l=k}^{K}P_l\cong \bigoplus_{l=k}^K F_l
\]
holds for $k=0,1,\ldots,K$ 
by the induction on $k$ from above. 
When $k=K$, it is OK. 
Suppose that we get the isomorphism for $k+1$. 
By adding the module $\bigoplus_{l=k+1}^{K}P_l$ 
to the sequence (\ref{FMPM}), 
we get the short exact sequence
\[
\begin{CD}
0 @>>> F_k @>>> M_k\oplus \bigoplus_{l=k}^{K}P_l @>>> 
M_{k+1}\oplus \bigoplus_{l=k+1}^{K}P_l @>>> 0. 
\end{CD}
\]
Since both $F_k$ and 
$M_{k+1}\oplus \bigoplus_{l=k+1}^{K}P_l\cong 
\bigoplus_{l=k+1}^K F_l$ 
are permutation modules, 
this sequence splits and we get the desired 
isomorphism 
$M_k\oplus \bigoplus_{l=k}^{K}P_l\cong \bigoplus_{l=k}^K F_l$ 
by Proposition~\ref{Prop:extpp}. 
Hence we obtain 
$M_0\oplus \bigoplus_{l=0}^{K}P_l\cong \bigoplus_{l=0}^K F_l$. 
Therefore $M_0$ is permutation projective. 
This completes the proof. 
\end{proof}

Combining Propositions~\ref{Prop:Ing3}, 
\ref{Prop:Ing1}, \ref{Prop:Ing2} and 
Theorem~\ref{Thm:TecMain}, 
we get Theorem~\ref{Thm:Main}.

\end{document}